# FRACTAL SUBSERIES OF THE HARMONIC SERIES

GÁBOR KORVIN[1]

ABSTRACT: We study the convergence of certain subseries of the harmonic series corresponding to increasing sequences of integers whose digits in a certain base are not uniformly distributed. We also discuss the case of *irregular sequences*, where the frequency distribution of some of the digits does not exist. Examples are given for irregular sequences where the corresponding harmonic subseries is convergent, or divergent, respectively.

## 1. INTRODUCTION

Kempner [1] observed that if we omit from the harmonic series $H := \sum_{n=1}^{\infty} \frac{1}{n}$ all terms corresponding to integers whose decimal representation contains at least once the digit 9, then the remaining "9-free" subseries $H_{10,9}$ (where 10 refers to the base, 9 is the missing digit) will be convergent, and bounded by 90. Hansberger's [2] *Mathematical Gems* gives a simple proof of this fact, and discusses the case of other missing digits. More recent papers [3-5] give better bounds, the best published result [5] being $20.2 \leq H_{10,9} \leq 28.3$.

In this paper we first we generalize Kempner's [1] observation (*Section 2*), and discuss the convergence properties of some *fractal subseries* of the harmonic series. Recall that a countable infinite set $\xi_i \in [0, \infty)$ of real numbers is called *fractal* with *mass fractal dimension* $\alpha$, $0 \leq \alpha < 1$ if their counting function $N(x) := \#(\xi_i \mid \xi_i \leq x)$ behaves for $x \to \infty$ as $N(x) \propto x^\alpha$ [6, 7]. In this spirit an increasing sequence $\{\nu_i\}$ of nonnegative integers will be called a *fractal sequence of dimension α* if their counting function asymptotically behaves as $N(n) := \#(v_i \mid v_i \leq n) \propto n^\alpha$ with some $\alpha$, $0 \leq \alpha < 1$.

---

[1] Box 1157, Earth Sciences Department, King Fahd University, Dhahran 31261, Saudi Arabia. Email: gabor@kfupm.edu.sa

*Section 3* is devoted to a more general case. The main result of this part (*Theorem 6* and its *Corollary 1*) studies the convergence of the subseries of the harmonic series corresponding to integers with *nonuniformly distributed digits* in their base-*r* expansion. The convergence will depend on the *Shannon entropy* of the digits' frequency distribution (if it exists!). Finally (in *Section 4*) examples will be constructed for increasing sequences of integers (called *irregular* sequences) for which the frequency distribution *does not exist* for certain digits. Surprisingly (see *Theorem 7*), for some irregular sequences the corresponding harmonic subseries is convergent, while in some other cases it is divergent. The concluding *Section 5* of the paper contains a simple conjecture, and related topics from *Number Theory* are indicated. Two technical *Lemmas* are relegated to the *Appendix*.

## 2. SIMPLE GENERALIZATIONS OF KEMPNER'S "9-FREE" SEQUENCES

*Theorem 1.* Let $v_{r,d} = \{v_i\}$ be the set of those integers which in their base-*r* expansion (*r>2*) do not contain a certain digit *d*, where *d* can be 0,1, …*r*-1. Then the counting function of $\{v_i\}$, $N_v(n) := \#(v_i \mid v_i \leq n)$ scales for $n \to \infty$ as $N_v(n) \propto n^\alpha$ with

$$\alpha = \frac{\log(r-1)}{\log r} + O\left(\frac{1}{\log n}\right).$$

*Remark 1.* If $r > 1$, the digits can be denoted by the symbols $\sigma_1 = 0, \sigma_2 = 1, \cdots, \sigma_{10} = 9$; and by arbitrary symbols for $\sigma_{11}, \cdots, \sigma_{r-1}$, as e.g. in the hexadecimal system one uses $0,1,\cdots,9, A, B, \cdots, F$.

*Proof of Theorem 1*: <u>Case a</u>) Let the missing digit be $d=0$. All integers $v$ such that $r^{k-1} \leq v \leq r^k$ have $k$ digits in base-$r$, that is the number of those $v$'s which do not contain the digit $d$ is $(r-1)^k$.

Letting $n = r^K$:

$$N_v(n) = \sum_{k=1}^{K} \#\left(\sum_{r^{k-1} \leq v < r^k} 1\right) = \sum_{k=1}^{K}(r-1)^k = -1 + \sum_{k=0}^{K}(r-1)^k$$

$$= -1 + \frac{(r-1)^{K+1} - 1}{r-2} \propto \frac{(r-1)^{K+1}}{r-2} \propto (r-1)^K = r^{\alpha K},$$

with $\alpha = \frac{\log(r-1)}{\log r} + O\left(\frac{1}{\log n}\right)$.

<u>Case b</u>) If the missing digit is $d \neq 0$, there are $(r-2)(r-1)^{k-1}$ possible $v$ values such that $r^{k-1} \leq v \leq r^k$, because the first digit cannot be 0. Thus, for $n = r^K$,

$$N_v(r^K) = \sum_{k=1}^{K}\#\left(\sum_{r^{k-1} \leq v < r^k} 1\right) = (r-2)\sum_{k=1}^{K-1}(r-1)^{k-1}$$

$$= (r-2)\sum_{k=0}^{K-2}(r-1)^k \propto (r-1)^K \propto r^{\alpha K},$$

again with $\alpha = \frac{\log(r-1)}{\log r} + O\left(\frac{1}{\log n}\right)$. ∎

*Definition* Let us call $\zeta_{r,d}(x) := \sum_{v \in \{v_{r,d}\}} \frac{1}{v^x}$ ($x$ real, $r>2$) the *generalized Euler zeta function* corresponding to the sequence $\{v_{r,d}\}$. We have the following

*Theorem 2.* $\zeta_{r,d}(x) := \sum_{v \in \{v_{r,d}\}} \frac{1}{v^x}$ *is convergent if* $x > \frac{\log(r-1)}{\log r}$ *and divergent if* $x \leq \frac{\log(r-1)}{\log r}$.

*Proof.* $\zeta_{r,d}(x) := \sum_{v \in \{v_{r,d}\}} \frac{1}{v^x} = \sum_{l=1}^{\infty} \sum_{r^{l-1} \le v < r^l} \frac{1}{v^x} < \sum_{l=1}^{\infty} \frac{(r-1)^l}{r^{(l-1)x}} = (r-1) \cdot \sum_{l=0}^{\infty} \frac{(r-1)^l}{r^{xl}}$ for *d=0* and

$$\zeta_{r,d}(x) := \sum_{v \in \{v_{r,d}\}} \frac{1}{v^x} = \sum_{l=1}^{\infty} \sum_{r^{l-1} \le v < r^l} \frac{1}{v^x} < \sum_{l=1}^{\infty} \frac{(r-2)(r-1)^{l-1}}{r^{xl}} < \sum_{l=0}^{\infty} \frac{(r-1)^l}{r^{xl}}$$ for *d≠0*.

In both cases if $x > \frac{\log(r-1)}{\log r}$ then $0 < \frac{r-1}{r^x} < 1$ and $\zeta_{r,d}(x) < \infty$. Suppose now that $x \le \frac{\log(r-1)}{\log r}$.

Then

$$\zeta_{r,d}(x) = \sum_{l=1}^{\infty} \sum_{r^{l-1} \le v < r^l} \frac{1}{v^x} > \sum_{l=1}^{\infty} \left( \sum_{r^{l-1} \le v < r^l} \frac{1}{r^{lx}} \right)$$

$$= \begin{cases} \sum_{l=1}^{\infty} \frac{(r-1)^l}{r^{xl}} & \text{if } d = 0 \\ \sum_{l=1}^{\infty} \frac{(r-2)(r-1)^{l-1}}{r^x \cdot r^{(l-1)x}} = \frac{r-2}{r^x} \cdot \sum_{l=0}^{\infty} \frac{(r-1)^l}{r^{lx}} & \text{if } d \ne 0, \end{cases}$$

and in both cases $\zeta_{r,d}(x)$ is divergent because $\frac{r-1}{r^x} \ge 1$. ∎

*Remark 2.* In classical Analysis ([8, vol. I, pp. 25-26] the exponent $\frac{\log(r-1)}{\log r}$ figuring in Theorem 2 is called the *convergence exponent* of the sequence $\{v_{r,d}\}$. For an arbitrary nondecreasing positive sequence $0 \le r_1 \le r_2 \le \cdots, \lim_{i \to \infty} r_i = \infty$, the convergence exponent is defined as the number λ for which $\sum_i r_i^{-\sigma} < \infty$ for $\sigma > \lambda$ and $\sum_i r_i^{-\sigma} = \infty$ if $\sigma < \lambda$. For $\sigma = \lambda$ the series may converge or diverge. It is proved in [8] that $\lambda = \limsup_{m \to \infty} \frac{\log m}{\log r_m}$.

*Theorem 2* is special case of the following more general Theorem.

*Theorem 3.* Let $\{v\} = \{v_1 < v_2 < \cdots\}$ be a monotone increasing sequence of integers with a counting function satisfying $N(n) \propto n^\alpha$, $\alpha \in (0,1]$ if $n > R$, where $R$ is a sufficiently large number. Then the generalized zeta function corresponding to the sequence $\{v\}$, $\varsigma_{\{v\}}(x) := \sum_i \frac{1}{v_i^x}$, is convergent if $x > \alpha$ and divergent if $x \leq \alpha$.

*Proof.* First we prove that if $x > \alpha$ then $\varsigma_{\{v\}}(x)$ is convergent. Indeed,

$$\varsigma_{\{v\}}(x) = \sum_{i=1}^\infty \frac{1}{v_i^x} = \sum_{l=1}^\infty \sum_{R^{l-1} \leq v_i < R^l} \frac{1}{v_i^x} < \sum_{l=1}^\infty \frac{N(R^l) - N(R^{l-1})}{(R^{l-1})^x} < \sum_{l=1}^\infty \frac{N(R^l)}{(R^{l-1})^x}$$

$$= const \cdot R^\alpha \sum_{l=0}^\infty \frac{R^{l\alpha}}{R^{lx}} < \infty \text{ if } R^{\alpha-x} < 1.$$

Next, let $x \leq \alpha$. We have, with appropriate positive constants $a$ and $b$:

$$\varsigma_{\{v\}}(x) = a + \sum_{l=2}^\infty \sum_{R^{l-1} \leq v < R^{l+1}} \frac{1}{v^x} > a + \sum_{l=1}^\infty \frac{N(R^{l-1}) - N(R^l)}{R^{lx}} > a + \sum_{l=1}^\infty \frac{N(R^{l-1}) - N(R^l)}{R^{(l+1)x}}$$

$$= a + \frac{b}{R^x} \sum_{l=1}^\infty \frac{R^{l\alpha}}{R^{lx}} = \infty \text{ if } \alpha \geq x.$$

∎

If *N(n)* is not a power function itself but asymptotically tends to such a function, we have a similar result:

*Theorem 4.* Let $\{v\} = \{v_1 < v_2 < \cdots\}$ be a monotone increasing sequence of integers with a counting function satisfying $N(n) \propto n^{\alpha[1+f(n)]}$, $\alpha \in (0,1]$, $\lim_{n \to \infty} f(n) = 0$. Then there exists an $\varepsilon_0$ such that for any $0 < \varepsilon < \varepsilon_0 < 1$ the generalized zeta function corresponding to the sequence, $\varsigma_{\{v\}}(x) := \sum_i \frac{1}{v_i^x}$, is convergent if $x > \alpha(1+\varepsilon)$ and divergent if $x \leq \alpha(1-\varepsilon)$.

*Proof.* Let $|f(n)| < \varepsilon_1 < 1$ for $n \geq R_1$ and $n^{\alpha(1-\varepsilon)} > 2$ for $n \geq R_2$ and $\varepsilon < \varepsilon_2 < 1$. Let $R = \max\{R_1, R_2\}$; $\varepsilon_0 = \min(\varepsilon_1, \varepsilon_2, 1)$. It is sufficient to study the convergence of the truncated series $\sum_{i>R} \frac{1}{v_i^x}$. First, let $x > \alpha(1+\varepsilon)$. Then

$$\sum_{i\geq R} \frac{1}{v_i^x} = \sum_{l=1}^{\infty} \sum_{R^l \leq v_i < R^{l+1}} \frac{1}{v_i^x} \leq \sum_{i=1}^{\infty} \frac{N_{\{v\}}(R^{l+1})}{R^{lx}} < const_1 \cdot \sum_{1}^{\infty} \frac{(R^{l+1})^{(1+\varepsilon)}}{R^{lx}} = const_2 \cdot \sum_{l=1}^{\infty} \frac{R^{l\alpha(1+\varepsilon)}}{R^{lx}} < \infty$$

if $\alpha(1+\varepsilon) - x < 0$. if $x \leq \alpha(1-\varepsilon)$ the scaling law $N(\lambda R^l) = \lambda^{\alpha(1+f[R^l])}$ implies

$$\sum_{l>L} \frac{1}{v_l^x} = \sum_{l>L} \sum_{R^l \leq v_i \leq R^{l+1}} \frac{1}{v_i^x} > \sum_{l=1}^{\infty} \frac{N(R^{l+1}) - N(R^l)}{R^{(l+1)x}} = \sum_{l=1}^{\infty} \frac{N(R^l)\left(R^{\alpha(1+f[R^l])} - 1\right)}{R^{(l+1)x}}$$

$$> \sum_{l=1}^{\infty} \frac{N(R^l)}{R^{(l+1)x}} \left(R^{\alpha(1-\varepsilon)} - 1\right) > const_3 \cdot \frac{1}{R^x} \sum_{l=1}^{\infty} \frac{R^{\alpha(1-\varepsilon)l}}{R^{lx}} = \infty \quad \text{if} \quad x \leq \alpha(1-\varepsilon)$$
∎

*Remark 3.* Theorems 3 and 4 can also be deduced from the following known *Theorem A* (of Krzyś, Powell and Šalát, [9, 10, 11]): *Let* $\{m_1 < m_2 < \cdots\}$ *be a sequence of integers with counting function N(n). Then* $\sum_{n=1}^{\infty} \frac{1}{m_n} < \infty$ *iff* $\sum_{n=1}^{\infty} \frac{N(n)}{n^2} < \infty$.

To derive *Theorem 3* from *Theorem A*, let sequence $\{v\} = \{v_1 < v_2 < \cdots\}$ be as in *Theorem 3*, and $N_{\{v\}}(n,x) := \#\{v_i^x \leq n\} = N_{\{v\}}\left(n^{\frac{1}{x}}\right) \propto n^{\frac{\alpha}{x}}$. By *Theorem A*, $\sum_{i=1}^{\infty} \frac{1}{v_i^x} < \infty$ iff

$\sum_{n=1}^{\infty} n^{\frac{\alpha}{x}-2} < \infty$ that is iff $x > \alpha$.

∎

To derive *Theorem 4* from *Theorem A*, let $N_{\{v\}}(n,x) := \#\{v_i^x \leq n\} = N_{\{v\}}\left(n^{\frac{1}{x}}\right) \propto n^{\frac{\alpha}{x}[1+f(n)]}$ and by *Theorem A* $\sum_{i=1}^{\infty} \frac{1}{v_i^x} < \infty$ iff $\sum_{n=1}^{\infty} n^{\frac{\alpha}{x}[1+f(n)]-2} < \infty$. Let $|f(n)| \leq \varepsilon < 1$ for $n \geq N_0$, then

$$\sum_{i=1}^{\infty} \frac{1}{v_i^x} = const + \sum_{n=N_0}^{\infty} \frac{1}{v_n^x} < \infty \quad \text{if} \quad \sum_{n \geq N_0} n^{\frac{\alpha}{x}[1+f(n)]-2} < \sum_{n \geq N_0} n^{\frac{\alpha}{x}[1+\varepsilon]-2} \quad \text{which is satisfied if}$$

$\alpha[1+\varepsilon] < x$. Conversely, if $x \leq \alpha[1-\varepsilon]$, we have

$$\sum_{n=1}^{\infty} n^{\frac{\alpha}{x}[1+f(n)]-2} > \sum_{n=1}^{N_0-1} n^{\frac{\alpha}{x}[1+f(n)]-2} + \sum_{n \geq N_0} n^{\frac{\alpha}{x}[1-\varepsilon]-2} = \infty. \quad \blacksquare$$

*Example* If $\{v\} = \{p_i\}$ is the sequence of prime numbers, then by the *Prime Number Theorem*

$$N_v(n) = \pi(n) = \int_2^n \frac{dx}{\log x} + O(x \cdot \exp(-a\sqrt{\log x})) \propto \frac{n}{\log n} \propto n^{1-\frac{\log \log n}{\log n}} \quad \text{for} \quad n \to \infty. \text{ By Theorem}$$

4, $\alpha = 1$, $f(n) = -\frac{\log \log n}{\log n} \to 0$, and $\zeta_{\{p_i\}}(x) = \sum_i \frac{1}{p_i^x} < \infty$ for $x > 1$ and divergent for $x < 1$.

Note that Theorem 4 cannot decide the convergence of the series for $x = 1$. Using the stronger *Theorem A*, $\sum_n \frac{\pi(n)}{n^2} \sim \sum_n \frac{1}{n \log n} = \infty$ and the well-known divergence of $\sum_n \frac{1}{p_n}$ follows.

*Theorem A* can also be used to prove the following Theorem.

*Theorem 5.* Let $\{v\} = \{v_1 < v_2 < \cdots\}$ be a sequence of non-negative integers with a counting function satisfying $N(n) \propto n^{\alpha[1+f(n)]}$, $\alpha \in (0,1]$. Suppose that for sufficiently large values of $n$ we have $\frac{n^{f(n)}}{(n+1)^{f(n+1)}} = 1 + \frac{\delta}{n} + O\left(\frac{1}{n^p}\right)$ with $p > 1$. Then $\sum_{i=1}^{\infty} \frac{1}{v_i^x} < \infty$ for $x > \alpha(1-\delta)$ and $\sum_{i=1}^{\infty} \frac{1}{v_i^x} = \infty$ for $x \leq \alpha(1-\delta)$.

*Proof.* By *Theorem A*, $\sum_{i=1}^{\infty} \frac{1}{v_i^x} < \infty$ iff $\sum_{n=1}^{\infty} \frac{N(n)}{n^2} \propto \sum_{n=1}^{\infty} n^{\frac{\alpha}{x}[1+f(n)]-2} < \infty$. According to the *Gauss criterium* ([12] Eq. 0.225, [13] p. 288), if for a sequence of positive numbers $\{u_k\}$ one has $\dfrac{u_k}{u_{k+1}} = 1 + \dfrac{q}{k} + \left|O\left(\dfrac{1}{k^\beta}\right)\right|$ with $\beta > 1$ then the series $\sum_{k=1}^{\infty} u_k$ is convergent for $q > 1$ and divergent for $q \leq 1$. Applying this to $\sum_{n=1}^{\infty} n^{\frac{\alpha}{x}[1+f(n)]-2}$ we have

$$\dfrac{n^{\frac{\alpha}{x}[1+f(n)]-2}}{(n+1)^{\frac{\alpha}{x}[1+f(n+1)]-2}} = \left(1+\dfrac{1}{n}\right)^{2-\frac{\alpha}{x}} \left\{\dfrac{n^{f(n)}}{(n+1)^{f(n+1)}}\right\}^{\frac{\alpha}{x}} = 1 + \left[2 + (\delta-1)\dfrac{\alpha}{x}\right] \cdot \dfrac{1}{n} + O\left(\dfrac{1}{n^p}\right). \blacksquare$$

3. SEQUENCES WITH NON-UNIFORMLY DISTRIBUTED DIGITS

In this part, we shall generalize *Theorems 1* and *2* for subseries of the harmonic series corresponding to integers with *non uniformly distributed digits* in their base-*r* expansion ($r \geq 2$). We denote by the symbols $\sigma_0, \sigma_1, \ldots, \sigma_{r-1}$ the digits $0, 1, \cdots$ used in the base-*r* number system. Let $\Lambda = (l_{ij})_{i=i_1, i_2, \cdots; j=0,1,\cdots,r-1}$ be a matrix of infinitely many rows and *r* columns, consisting of nonnegative integer elements. Assume that $\Lambda$ has the following properties:

(A) $\sum_{j=0}^{r-1} l_{ij} = i$ for $i = i_1, i_2, \cdots$ with $i_1 \leq i_2 \leq \cdots$; $\lim_{n\to\infty} i_n = \infty$,

(B) $\lim_{i\to\infty} \dfrac{l_{ij}}{i} := \lambda_j$ exists for all $j = 0, 1, \cdots, r-1$, and

(C) $\lambda_j \neq 0$ for $j = 0, 1, \cdots, r-1$.

Then, of course, we also have

(D) $\sum_{j=0}^{r-1} \lambda_j = 1$.

Denote by $V_\Lambda$ the set of all positive integers $v$ for which:

*(E)* If $r^{i-1} \leq v < r^i$ then the base-$r$ expansion of $v$ contains the digits $\sigma_0, \sigma_1, \ldots, \sigma_{r-1}$ exactly $l_{i0}, l_{i1}, \cdots, l_{i,r-1}$ times ($i = i_1, i_2, \cdots$), respectively;

*(F)* The most significant (i.e. left-most) digit of $v$ is not zero.

*Theorem 6.* $\zeta_{V_\Lambda}(x) := \sum\limits_{v \in V_\Lambda} \dfrac{1}{v^x}$ is convergent if $x > \dfrac{-\sum\limits_{j=0}^{r-1} \lambda_j \log \lambda_j}{\log r}$ and divergent if

$$x \leq \dfrac{-\sum\limits_{j=0}^{r-1} \lambda_j \log \lambda_j}{\log r}.$$

*Remark 4.* As the maximum of $-\sum\limits_{j=0}^{r-1} \lambda_j \log \lambda_j$, subject to the conditions that $\lambda_j \geq 0$ for $j = 0, 1, \cdots, r-1$ and $\sum\limits_{j=0}^{r-1} \lambda_j = 1$, is $\log r$, *Theorem 6* has the following Corollary:

*Corollary 1.* The $\Lambda$-subseries of the harmonic series, $H_{V_\Lambda} = \sum\limits_{v \in V_\Lambda} \dfrac{1}{v}$ is divergent if $\lambda_0 = \lambda_1 = \cdots = \lambda_{r-1} = \dfrac{1}{r}$ and it is convergent otherwise.

To prove Theorem 6 we shall need the following Lemma:

*Lemma 1.* Denote by $\Theta_r$ the set of those integers which in their base-$r$ expansion ($r \geq 2$) do not contain all symbols $\sigma_0 = 0, \sigma_1 = 1, \sigma_2, \ldots, \sigma_{r-1}$ at least once. Then the series $\sum\limits_{v \in \Theta_r} \dfrac{1}{v}$ is convergent.

*Proof of Lemma 1.* As in *Theorem 1*, we denote by $V_{r,\sigma} = \{v\}$ the set of those integers which in their base-$r$ expansion ($r \geq 2$) do not contain a certain digit $\sigma$, where $\sigma$ can be any of $\sigma_0 = 0, \sigma_1 = 1, \sigma_2, \ldots, \sigma_{r-1}$. Also, as in *Theorem 2*, we use the notation $\zeta_{\Theta_r}(x) := \sum_{v \in \Theta_r} \frac{1}{v^x}$. Then we have $\sum_{v \in \Theta} \frac{1}{v} < \sum_{k=0}^{r-1} \sum_{v \in V_{r,\sigma_k}} \frac{1}{v} = \sum_{k=0}^{r-1} \zeta_{r,\sigma_k}(1) < \infty$ because by

*Theorem 2* for all $k = 0, 1, \cdots, r-1$ $\zeta_{r,\sigma_k}(1) < \infty$. ∎

*Proof of Theorem 6.* Write $v_\Lambda = (v_\Lambda \cap \Theta_r) \cup (v_\Lambda \setminus \Theta) \subseteq \Theta_r \cup v_\Lambda^*$ where the set $v_\Lambda^*$ consists of those numbers from $v_\Lambda$ which contain in their base-$r$ expansion ($r \geq 2$) all possible digits at least once. Observe, that all these numbers are necessarily at least $r$ digits long, that is if $v \in v_\Lambda^*$ then $v \geq r^{r-1}$. Because of *Lemma 1* it is sufficient to study the convergence of the series corresponding to the subsequence $v_\Lambda^*$. First we derive an upper bound for $N(r,i) := \#\{v \mid v \in v_\Lambda^*, r^{i-1} < v < r^i\}$ where $i > r-1$. All counted numbers contain $i$ digits. By the definition of $v_\Lambda^*$ we have $l_{i0} \geq 1, l_{i1} \geq 1, \cdots, l_{i,r-1} \geq 1$. We will get an upper bound $\overline{N}(r,i) > N(r,i)$ if we also count such numbers $v \in v_\Lambda^*$ whose first digit is 0, that is if we write $N(r,i) < \overline{N}(r,i) = \frac{i!}{l_{i0}! l_{i1}! \cdots l_{i,r-1}!}$. Because of conditions B and C, there exists a sufficiently large $I$ such that for $i > I$ even $n = \min\{l_{i0}, \cdots, l_{i,r-1}\}$ will be as large that we can apply Stirling's formula $n! \sim \left(\frac{n}{e}\right)^n \sqrt{2\pi n}$. The formula gives

$$\log N(r,i) < \log i! - \sum_{j=0}^{r-1} \log l_{ij}! \approx \left[i \log i - i + \frac{1}{2} \log 2\pi i\right] - \left[\sum_{j=0}^{r-1} l_{ij} \log l_{ij} - \sum_{j=0}^{r-1} l_{ij} + \frac{1}{2} \sum_{j=0}^{r-1} \log 2\pi l_{ij}\right]$$

$$= -i \sum_{j=0}^{r-1} \frac{l_{ij}}{i} \log \frac{l_{ij}}{i} + \frac{1}{2} \log 2\pi i - \frac{1}{2} \sum_{j=0}^{r-1} \log 2\pi l_{ij} < -i \sum_{j=0}^{r-1} \frac{l_{ij}}{i} \log \frac{l_{ij}}{i} + \frac{1}{2} \log 2\pi i,$$

that is

$$N(r,i) < \overline{N}(r,i) < \sqrt{2\pi i} \exp\left[-i \sum_{j=0}^{r-1} \frac{l_{ij}}{i} \log \frac{l_{ij}}{i}\right] = \sqrt{2\pi i} \exp\left[i \cdot E\left(\frac{l_{i0}}{i}, \cdots, \frac{l_{i,r-1}}{i}\right)\right] \quad (1)$$

where $E\left(\dfrac{l_{i0}}{i},\cdots,\dfrac{l_{i,r-1}}{i}\right) = -\sum_{j=0}^{r-1}\dfrac{l_{ij}}{i}\log\dfrac{l_{ij}}{i}$ is the *Shannon entropy* of the probability distribution $\left(\dfrac{l_{i0}}{i},\cdots,\dfrac{l_{i,r-1}}{i}\right)$. By Condition *B*, and the continuity of the Shannon entropy, for an arbitrarily small positive ε there exists an $i_\varepsilon$ such that for $i > i_\varepsilon$ one has $E\left(\dfrac{l_{i0}}{i},\cdots,\dfrac{l_{i,r-1}}{i}\right) \leq (1+\varepsilon)E(\lambda_0,\lambda_1,\cdots,\lambda_{r-1})$ where $\lambda_j = \lim_{i\to\infty}\dfrac{l_{ij}}{i}$. For such an *i* we have

$$N(r,i) < \sqrt{2\pi i}\,\exp[i\cdot E(\lambda_0,\cdots,\lambda_{r-1})]\ . \qquad (2)$$

Introduce the notations $i^* = \max\{i_\varepsilon, r\}$, $E(\lambda_0,\lambda_1,\cdots,\lambda_{r-1}) = E$, $\sqrt{2\pi} = B$. Then

$$\sum_{v\in V_\Lambda^*}\dfrac{1}{v} < A + \sum_{i=i^*}^{\infty}\sum_{\substack{v\in V_\Lambda^*\\ r^{i-1}\leq v < r^i}}\dfrac{1}{v} < A + \sum_{i=i^*}^{\infty}\dfrac{N(r,i)}{r^{i-1}} < A + Br\sum_{i=i^*}^{\infty}\dfrac{\sqrt{i}}{r^i}\exp(iE)\ . \qquad (3)$$

As the maximum of $E(\lambda_0,\lambda_1,\cdots,\lambda_{r-1})$, subject to the conditions that $\lambda_j \geq 0$ for $j = 0,1,\cdots,r-1$ and $\sum_{j=0}^{r-1}\lambda_j = 1$, is $E_{\max} = E(\dfrac{1}{r},\cdots,\dfrac{1}{r}) = \log r$, Eq. (3) can be written as

$$\sum_{v\in V_\Lambda^*}\dfrac{1}{v} < A + Br\sum_{i=i^*}^{\infty}\dfrac{\sqrt{i}}{r^i}\exp(iE) = A + Br\sum_{i=i^*}^{\infty}\sqrt{i}\,\exp i(E - E_{\max}). \qquad (4)$$

Consider first the case when $E - E_{\max} < 0$. Then the series $\sum_{i=i^*}^{\infty}\sqrt{i}\,\exp i(E - E_{\max})$ is convergent by the Cauchy ratio test, because for sufficiently large values of *i* we have

$$\sqrt{\dfrac{i+1}{i}}\exp(E - E_{\max}) < 1.$$

An argument similar to that which has lead to inequality (4) gives, with different constants $A'$ and $B'$ and a different starting summation index $i^{\bullet\prime}$ that

$$\zeta_{V_\Lambda}(x) := \sum_{v \in V_\Lambda} \frac{1}{v^x} < A' + B r^x \sum_{i=i^*}^{\infty} \sqrt{i} \exp i(E - xE_{max}) \ . \tag{5}$$

By Cauchy's ratio test $\zeta_{V_\Lambda}(x) := \sum_{v \in V_\Lambda} \frac{1}{v^x} < \infty$ if $E - xE_{max} = E - x\log r < 0$ that is if

$$x > \frac{E}{\log r} = \frac{-\sum_{j=0}^{r-1} \lambda_j \log \lambda_j}{\log r} \ . \tag{6}$$

The case $E = E_{max}$ must be differently handled. In this case the limiting distribution of the $\Lambda$-matrix elements is $\lambda_0 = \lambda_1 = \cdots = \lambda_{r-1} = \frac{1}{r}$. By condition B and the continuity of the Shannon entropy, for an arbitrarily small $\delta > 0$ there exists an index $i^{**}$ such that

$$l_{ij} > 0 \quad for \quad i \geq i^{**}, j = 0,1,\cdots,r-1 \tag{7.a}$$

$$E\left(\frac{l_{i0}}{i}, \cdots, \frac{l_{i,r-1}}{i}\right) \geq (1-\delta) E\left(\frac{1}{r}, \cdots, \frac{1}{r}\right) = (1-\delta) \log r \ . \tag{7.b}$$

Denote the corresponding sequence, for $i \geq i^{**}$, by $v_\Lambda^*$. We need a lower bound $\underline{N}(r,i)$ for

$N(r,i) = \#\{v \mid v \in v_\Lambda^*, r^{i-1} < v < r^i\}$, $i \geq i^{**}$. Because of $\lambda_0 = \lambda_1 = \cdots = \lambda_{r-1} = \frac{1}{r}$, a fraction $\frac{r-1}{r}$ of the $\frac{i!}{l_{i0}! l_{i1}! \cdots l_{i,r-1}!}$ formally possible numbers do not contain 0 in their most significant digit. To find $\underline{N}(r,i)$ we use Eq. (1) with $(1-\delta)E_{max} = (1-\delta)\log r$ instead of $E\left(\frac{l_{i0}}{i}, \cdots, \frac{l_{i,r-1}}{i}\right)$, which gives

$N(r,i) > \underline{N}(r,i) = \frac{r-1}{r} \sqrt{2\pi} (1-\delta) \exp(i \log r) > (1-\delta) r^i$ for $r \geq 2$. Consequently,

$$\sum_{v \in V_{\Lambda^*}} \frac{1}{v} > \sum_{i=i^{**}}^{\infty} \sum_{\substack{v \in V_{\Lambda^*} \\ r^{i-1} \leq v < r^i}} \frac{1}{r^{i-1}} > \sum_{i=i^{**}}^{\infty} \frac{N(r,i)}{r^i} > (1-\delta) \sum_{i=i^{**}}^{\infty} \frac{\sqrt{i}\, r^i}{r^i} = \infty.$$

The same argument as before leads, with an appropriate constant $B''>0$ and starting index,

to $\zeta_{V_{\Lambda^*}}(x) > B'' \sum_{i=(i^{**})'}^{\infty} \sqrt{i} \exp[\log r^{i(1-x)}] = B'' \sum_{i=(i^{**})'}^{\infty} \sqrt{i}\, r^{i(1-x)}$ and by Cauchy's ratio ratio test

$\zeta_{V_{\Lambda}^*} = \infty$ if $x \leq 1$.

Finally, we show that $\zeta_{V_{\Lambda}}(x) = \infty$ if $x \leq \frac{E(\lambda_0, \lambda_1, \cdots, \lambda_r)}{\log r}$. In this case a fraction $1-\lambda_0$ of the numbers $v \in V_{\Lambda}$ do not contain 0 as their most significant digit. With some appropriate constant $B''' > (1-\lambda_0)\sqrt{2\pi}$ and for a sufficiently large index $i^{***}$ the previous argument gives $\zeta_{V_{\Lambda}}(x) > B''' \sum_{i=i^{***}}^{\infty} \sqrt{i} \exp[i(E - xE_{\max})]$. By the ratio test, the RHS is divergent if

$$\frac{E}{E_{\max}} \geq x. \qquad \blacksquare$$

## 4. IRREGULAR SEQUENCES

We emphasize that condition $B$ in the definition of $\Lambda$ is not trivial, because of the following Theorem:

*Theorem 7. There exist increasing sequences of integers (in base r) for which the limit frequency* $\lim_{i \to \infty} \frac{l_{ij}}{i} := \lambda_j$ *does not exist for some* $j \in \{1, 2, \cdots, r-1\}$. *We call them irregular sequences. If r>2, there exist such irregular sequences* $V_{\Lambda}^{conv}$ *with the further property that* $\sum_{v \in V_{\Lambda}^*} \frac{1}{v} < \infty$, *and there exist irregular sequences* $V_{\Lambda}^{div}$ *for which* $\sum_{v \in V_{\Lambda}^*} \frac{1}{v} = \infty$.

*Remark 5.* The term "irregular" in the Theorem has been borrowed from papers of Barreira *et al.* [14, 15] who called a real number $\alpha \in (0,1)$ *irregular* or *atypical* if the frequency distribution of digits in their decimal (or in any other base *r*) expansion does not exist.

*Proof.* To construct a simple irregular sequence $v_\Lambda^{conv}$ in base *r*, select any of the *r* possible digits, say $\sigma_{j_0} = \sigma$ and replace all other digits by the symbol $*$. Consider the sequence

$$\begin{cases} v_0 = 1 \\ v_1 = 1\sigma \\ v_2 = 1\sigma * \\ v_3 = 1\sigma ** \\ v_4 = 1\sigma **\sigma \\ v_5 = 1\sigma **\sigma\sigma \\ v_6 = 1\sigma **\sigma\sigma\sigma \\ v_7 = 1\sigma **\sigma\sigma\sigma\sigma \\ v_8 = 1\sigma **\sigma\sigma\sigma\sigma * \\ v_9 = 1\sigma **\sigma\sigma\sigma\sigma ** \\ \vdots \\ v_{15} = 1\sigma **\sigma\sigma\sigma\sigma ******** \\ v_{16} = 1\sigma **\sigma\sigma\sigma\sigma ********\sigma \\ \vdots \\ v_{31} = 1\underbrace{\sigma}_{1\times}\underbrace{*\,*}_{2\times}\underbrace{\sigma\sigma\sigma\sigma}_{4\times}\underbrace{********}_{8\times}\underbrace{\sigma\cdots\sigma}_{16\times} \\ \vdots \end{cases}$$

It is easy to check that $\lim\limits_{i\to\infty}\dfrac{l_{i\sigma}}{i} := \lambda_\sigma$ (where $l_{i\sigma}$ is the number of the $\sigma$-digits in the number $v_i$) does not exist, because the limit is:

$$\lim_{i\to\infty}\frac{l_{i\sigma}}{i} = \frac{1}{3} \text{ if } i \to \infty \text{ as } i = 1, 7, 31, \cdots, 2^{2k-1} - 1, \cdots \ (k = 1, 2, 3, \cdots)$$

$$\lim_{i\to\infty}\frac{l_{i\sigma}}{i}=\frac{2}{3} \text{ if } i\to\infty \text{ as } i=3,15,63,\cdots,2^{2k}-1,\cdots \quad (k=1,2,3,\cdots).$$

More generally, with the notation $\sigma^n = \underbrace{\sigma\sigma\cdots\sigma}_{n-times}$, $n \geq 1$, we can construct the sequence

$$\{v\} = v_\Lambda^{(k)} = \begin{cases} v_0 = \gamma \\ v_1 = \gamma\sigma \\ v_2 = \gamma\sigma * \\ v_3 = \gamma\sigma *^2 \\ \vdots \\ v_{k+1} = \gamma\sigma *^k \\ v_{k+2} = \gamma\sigma *^k \sigma \\ \vdots \\ v_{k^2+k+1} = \gamma\sigma *^k \sigma^{k^2} \\ \vdots \\ v_{k^3+k^2+k+1} = \gamma\sigma *^k \sigma^{k^2} *^{k^3} \\ \vdots \end{cases}$$

where $* \neq \sigma$, $\gamma \neq \sigma, 0$ can be arbitrary digits. It is easy to check (see *Lemma 2* in the Appendix) that

$$\liminf_{i\to\infty}\frac{l_{i\sigma}}{i}=\frac{1}{k+1}; \limsup_{i\to\infty}\frac{l_{i\sigma}}{i}=\frac{k}{k+1}, \tag{8}$$

that is the sequence $\left\{\frac{l_{i\sigma}}{i}\right\}$ is not convergent if $k > 1$. There are at least two convergent subsequences of $\left\{\frac{l_{i\sigma}}{i}\right\}$ which tend to different limits, namely to $\frac{1}{k+1}$ and $\frac{k}{k+1}$. These are obtained when $i \to \infty$ through the numbers $i = \dfrac{k^{2m+2}-1}{k-1}$, or through the numbers $i = \dfrac{k^{2m+4}-1}{k-1}$, $(m=1,2,3,\cdots)$, respectively.

If $r>2$ and we substitute for all $*$-digits the value $\gamma$ $(\gamma \neq \sigma,\ \gamma \neq 0)$, we get an *irregular sequence* $V_\Lambda^{k,conv}$. (The convergence of the corresponding series $\sum_{v \in V_\Lambda^{k,conv}} \dfrac{1}{v}$ follows from *Lemma 1*, because $V_\Lambda^{k,conv} \subseteq \Theta_r$.)

Next we construct a sequence $V_\Lambda^{k,div}$ for which (when represented in base $r$, $r>2$) $\lim_{i \to \infty} \dfrac{l_{ij}}{i} := \lambda_j$ does not exist for some $j \in \{1,2,\cdots,r-1\}$, and the sum $\sum_{v \in V_\Lambda^{k,div}} \dfrac{1}{v}$ is divergent.

Starting out from the previous sequence $\{v\} = V_\Lambda^{k,conv}$ we can include further numbers to it if we select the $*$-digits less restrictively, for example by allowing each of them to be any symbol except $\sigma$. We sort the terms of the sequence in increasing order so that $v_1 < v_2 < \cdots$, define $i_n$ as the length of (that is number of digits in) the number $v_n$, and $l_{i_n,\sigma}$ as the number of $\sigma$-digits in $v_n$. Then the indexes $\{i_n\}$ are nondecreasing. Denote the corresponding sequence by $V_\Lambda^{(k)}$. By *Lemma 3* (in the Appendix) the *lim inf* and *lim sup* of the sequence of real numbers $\left\{\dfrac{l_{i_n,\sigma}}{i_n}, n=1,2,\cdots\right\}$ are the same as of the sequence $\left\{\dfrac{l_{i\sigma}}{i}\right\}$.

The corresponding extended series $\sum_{v \in V_\Lambda^{(k)}} \dfrac{1}{v}$ however is still convergent. To see this, estimate its counting function $M(N,k) := \#\{v \mid v \in V_\Lambda^{(k)}, v \leq N\}$. If $r>2$, there are $(r-2)$ ways to select the first digit $\gamma$ of the numbers in the $V_\Lambda^{(k)}$ sequence, and $(r-1)$ ways to select each subsequent $*$-digit. By Eq. (8), and *Lemmas 2* and *3*, for almost all $v \in V_\Lambda^{(k)}$ consisting of $M$ ($M = r^I$ with $I \gg 1$) digits the number of $*$-digits in $v$ is less than or equal to $\dfrac{Mk}{k+1}$. Consequently, the counting function $M(N,k) = \#\{v \mid v \in V_\Lambda^{(k)}, v \leq N\}$ satisfies

$$M(N,k) = \#\{v \mid v \in V_\Lambda^{(k)}, v \leq N\} = \sum_{i=1}^{I}\{v \mid v \in V_\Lambda^{(k)}, r^{i-1} \leq v < r^i\} \leq (r-2)\sum_{i=1}^{I}\left[(r-1)^{\frac{k}{k+1}}\right]^i$$

$$< (r-2)\frac{(r-1)^{I+1} - 1}{r-2} < (r-1)^{I+1} = (r-1)(r-1)^{\frac{\log N}{\log r}} = (r-1) N^{\frac{\log(r-1)}{\log r}} = (r-1) N^{1-\varepsilon}$$

where $\varepsilon = \dfrac{1}{r \log r}$, $0 < \varepsilon < 1$. As the series $\sum_{N=1}^{\infty} \dfrac{N^{1-\varepsilon}}{N^2} = \sum_{N=1}^{\infty} \dfrac{1}{N^{1+\varepsilon}}$ is convergent, by *Theorem A* (of Krzyś, Powell and Šalát, [9, 10, 11]) $\sum_{v \in V_\Lambda^{(k)}} \dfrac{1}{v} < \infty$.

To get a divergent irregular series we should increase the number of terms much more aggressively. Distribute the approximately $\dfrac{M}{k+1}$ $\sigma$-digits among the $M$ digits in all possible ways, there are $\begin{pmatrix} M \\ \dfrac{M}{k+1} \end{pmatrix}$ different ways to do this; and select the remaining $\dfrac{Mk}{k+1}$ $*$-digits arbitrarily, which can be done in $\left(\dfrac{Mk}{k+1}\right)^{r-1}$ ways. We obtain

$$\begin{pmatrix} M \\ \dfrac{M}{k+1} \end{pmatrix}\left(\dfrac{Mk}{k+1}\right)^{r-1} \sim k^M M^{r-1}\left(1 - \dfrac{r-1}{k}\right) \sim k^M M^{r-1}$$ different $M$-digit numbers if $N \gg k \gg 1$. Select $N$ as $N = r^I$. We have

$$M(N,k) = \sum_{i=1}^{I}\left\{v \mid v \in V_\Lambda^{(k)}, r^{i-1} \leq v < r^i\right\} \sim \sum_{i=1}^{I} k^i i^{r-1} < I^{r-1}\sum_{i=1}^{I} k^i = I^{r-1}\frac{k^{I+1}-1}{k-1}$$

$$\sim I^{r-1}k^I \sim (\log N)^{r-1}\left[k^{1/\log r}\right]^{\log N},$$

that is, by *Theorem A* (of Krzyś, Powell and Šalát, [9, 10, 11]), the convergence of the corresponding series $\sum \frac{1}{v}$ is equivalent to that of $\sum_{N=1}^{\infty}(\log N)^{r-1}\left[k^{1/\log r}\right]^{\log N} N^{-2}$.

Checking of the *Gauss criterium* ([12, 13]) reduces to the study of

$$\left(\frac{\log N}{\log(N+1)}\right)^{r-1} \alpha^{\log N - \log(N+1)}\left(\frac{N+1}{N}\right)^2 \text{ where } \alpha = k^{1/\log r}.$$

But $\left(\frac{\log N}{\log(N+1)}\right)^{r-1} \approx 1 - \frac{r-1}{N\log N}$, $\alpha^{\log N - \log(N+1)} \approx 1 - \frac{\log \alpha}{N}$, $\left(\frac{N+1}{N}\right)^2 \approx 1 + \frac{2}{N}$, that is

$$\left(\frac{\log N}{\log(N+1)}\right)^{r-1} \alpha^{\log N - \log(N+1)}\left(\frac{N+1}{N}\right)^2$$

$$\approx \left[1 - \frac{r-1}{N\log N}\right]\left[1 - \frac{\log \alpha}{N}\right]\left[1 + \frac{2}{N}\right] \approx 1 + \frac{1}{N}\left[2 - \frac{r-1}{\log N} - \log \alpha\right] + O(N^{-2})$$

$$= 1 + \frac{1}{N}\left[2 - \frac{r-1}{\log N} - \frac{\log k}{\log r}\right] + O(N^{-2})$$

For $k < r^2$, and a sufficiently large $N$, $q = \left[2 - \frac{r-1}{\log N} - \frac{\log k}{\log r}\right] > 1$ and by the Gauss criterium $\sum_{N=1}^{\infty}\frac{M(N,k)}{N^2} < \infty$. By *Theorem A* we also have $\sum \frac{1}{v} < \infty$. For $k > r^2$ however we have $\sum_{N=1}^{\infty}\frac{M(N,k)}{N^2} > \infty$ so that $\sum_v \frac{1}{v} = \infty$. ∎

## 5. A CONJECTURE AND FINAL REMARKS

The use of the Stirling formula in the proof of *Theorem 6* necessitated the assumption that in the limiting probability distribution of $\Lambda$, $\{\lambda_0, \cdots, \lambda_{r-1}\}$, all $\lambda_j$ are positive (condition *C*). It seems likely that *Theorem 6 remains valid even if some of the $\lambda_i$ are zero.*

We are aware of the formal similarities between the questions discussed in this paper and the probabilistic theory of the distribution of the decimal digits of irrational numbers (Borel's normal numbers and related topics, [16-19]), of the dimensional theory of continued fractions [20], or the Hausdorff dimension of the random Cantor set [21, 22]. Discussions of this similarity however would be outside the scope of this elementary paper.

*Acknowledgments.* The problem was suggested by a question of Prof. Klavdia Olechko (*UNAM*, Mexico City) about fractals and singularities, and this paper is dedicated to her. The happy and peaceful atmosphere, and excellent research facilities, of the *King Fahd University* where this research has been done, are gratefully acknowledged.

REFERENCES

[1] A.J. Kempner, "A curious convergent series", *American Math. Monthly* 1914: 48-50

[2] W. Stadje, "Konvergenz von Teilen der harmonischen Reihe", *Elem. Math.* **46**(1991) No. 2: 51-54

[3] R. Hansberger, *Mathematical Gems,* vol. 2. Math. Assoc. America, Wash., DC, 1976, pp. 98-103.

[4] A.D. Wadhwa, "Some convergent subseries of the harmonic series", *American Math. Monthly* **85**(1978) No 8: 661-665.

[5] A.D. Wadhwa, "An interesting subseries of the harmonic series", *American Math. Monthly* **82**(1975) No. 9: 931-933.

[6] G. Korvin, *Fractal Models in the Earth Sciences,* Elsevier, Amsterdam, 1992.

[7] J. Feder, *Fractals,* Plenum Press, New York, 1988.

[8] G. Pólya and G. Szegö, *Problems and Theorems in Analysis,* Springer Verlag, Berlin-Heidelberg-New York, 1978, vol. I, pp. 25-26.

[9] J. Krzyś, "Oliver's theorem and its generalizations", *Prace Mat.,* **2**(1956): 159-164.

[10] T. Šalát, *Infinite Series,* Academia, Prague, 1974, p. 101.

[11] B.J. Powell and T. Šalát, "Convergence of subseries of the harmonic series and asymptotic densities of sets of integers", *Publ. Inst. Math. (Beograd)* **50**(64)1991: 60-70.

[12] I.S. Gradshtein and I.M. Ryzhik, *Tables of Integrals, Sums, Series and Products*. Gos. Izdatel'stvo Fiziko-Matematicheskoi Literatury, Moscow, 1963 (In Russian).

[13] K. Knopp, *Theory and Applications of Infinite Series,* Hafner Publ. Co., New York, 1971 p. 288.

[14] L. Barreira, B. Saussol, and J. Schmeling, "Distribution of frequencies of digits via multifractal analysis", *J. Number Theory* **97**(2002): 413-442.

[ 15 ] L. Barreira and J. Schmeling, "Sets of "non-typical" points have full topological entropy and full Hausdorff dimension", *Isr. J. Math.* **116**(2000): 29-70.

[16] E. Borel, "Sur les probabilités dénombrables et leurs applications arithmétiques", *Rend. Circ. Mat. Palermo*, **26**(1909): 247-271.

[17] I. Niven, *Irrational Numbers,* The Carus Mathematical Monographs, John Wiley & Sons, NY, 2nd edn., 1963.

[18] P. Billingsley, *Ergodic Theory and Information,* Robert E. Krieger Publ. Co., Huntington, NY, 1978.

[19] H. G. Eggleston, "The fractional dimension of a set defined by decimal properties", *Quart. J. Math. Oxford Ser.,* **20**(1949): 31-36.

[20] I.J. Good, "The fractional dimensional theory of continued fractions", *Proc. Cambr. Phil. Soc.,* **37**(1941): 199-228.

[21] F. Hausdorff, "Dimension und äusseres Mass", *Math. Annalen,* **79**(1919): 157-179.

[22] K. Falconer, *Fractal Geometry. Mathematical Foundations and Applications.* John Wiley and Sons, Chichester, UK, 1985.

APPENDIX

*Lemma 2. Consider the sequence of increasing integers in base r (r>2)*

$$v_\Lambda^{(k)} = \begin{cases} v_0 = \gamma \\ v_1 = \gamma\sigma \\ v_2 = \gamma\sigma * \\ v_3 = \gamma\sigma *^2 \\ \vdots \\ v_{k+1} = \gamma\sigma *^k \\ v_{k+2} = \gamma\sigma *^k \sigma \\ \vdots \\ v_{k^2+k+1} = \gamma\sigma *^k \sigma^{k^2} \\ \vdots \\ v_{k^3+k^2+k+1} = \gamma\sigma *^k \sigma^{k^2} *^{k^3} \\ \vdots \end{cases}$$

*where $k \geq 2$, the digit $\gamma$ is different from $\sigma$ and 0, the digit $*$ is not the same as $\sigma$, and we use the notation $\sigma^n = \underbrace{\sigma\sigma\cdots\sigma}_{n-times}$, $n \geq 1$. Denote by $l_{i\sigma}$ the number of $\sigma$ digits in $v_i$ $(i = 1, 2, \cdots)$. Then*

$$\liminf_{i \to \infty} \frac{l_{i\sigma}}{i+1} = \frac{1}{k+1}; \limsup_{i \to \infty} \frac{l_{i\sigma}}{i+1} = \frac{k}{k+1} \quad . \tag{A.1}$$

*Proof.* To prove the first part of (A.1) we show that for any $\varepsilon > 0$ there exists an integer $N = N(\varepsilon, k)$ such that for all integers $n \geq N$ the following two conditions are satisfied:

(i) $\dfrac{l_{n\sigma}}{n+1} > \dfrac{1}{k+1} - \varepsilon$ ;

(ii) there exists an integer $p$, $p \geq n$ for which $\dfrac{l_{p\sigma}}{p+1} < \dfrac{1}{k+1} + \varepsilon.$

Select $N$ as $N = 1 + k + k^2 + \ldots + k^{2m+1} = \dfrac{k^{2m+2} - 1}{k - 1}$ for a sufficiently large $m$, where $m = m(k, \varepsilon)$ will be specified at the end of the proof, see Eqs. (A.7) and (A.8). Let $n > N$.

There exists a unique integer $m' \geq m$ for which $\dfrac{k^{2m'+2}-1}{k-1} < n \leq \dfrac{k^{2m'+4}}{k-1}$. We distinguish between two cases:

$$a) \quad \frac{k^{2m'+2}-1}{k-1} < n \leq \frac{k^{2m'+3}-1}{k-1}$$

$$b) \quad \frac{k^{2m'+3}}{k-1} < n \leq \frac{k^{2m'+4}-1}{k-1}$$

In case *a)*, because of the construction of sequence $v_\Lambda^{(k)}$, the total number of digits in $v_n$ is less than $\dfrac{k^{2m'+3}-1}{k-1}+1$, the number of $\sigma$-digits is greater than $\dfrac{k^{2m'+4}-1}{k^2-1}$, that is

$$\frac{l_{n\sigma}}{n+1} > \frac{k^{2m'+4}-1}{k^2-1} \cdot \left(\frac{k^{2m'+3}-1}{k-1}+1\right)^{-1} = \frac{1}{k+1} \cdot \frac{k^{2m'+4}-1}{k^{2m'+3}-1}\left(1+\frac{k-1}{k^{2m'+4}}\right)^{-1}. \tag{A.2}$$

As the right-hand-side of Eq. (A.2) tends to $\dfrac{k}{k+1}$ as $m' \to \infty$, for sufficiently large values of $m'$ we shall have

$$\frac{l_{n\sigma}}{n+1} > \frac{k}{k+1} - \varepsilon > \frac{1}{k+1} - \varepsilon. \tag{A.3}$$

In case *b)*, because of the construction of sequence $v_\Lambda^{(k)}$, the total number of digits in $v_n$ is less than $\dfrac{k^{2m'+4}-1}{k-1}+1$, the number of $\sigma$-digits is identically $\dfrac{k^{2m'+4}-1}{k^2-1}$ for all $n$ in this range, that is for sufficiently large values of $m'$

$$\frac{l_{n\sigma}}{n+1} > \frac{k^{2m'+4}-1}{k^2-1}\cdot\left(\frac{k^{2m'+4}-1}{k-1}+1\right)^{-1} = \frac{k^{2m'+4}-1}{k^2-1}\cdot\frac{k-1}{k^{2m'+4}-1}\left(1+\frac{k-1}{k^{2m'+4}}\right)^{-1}$$
$$> \frac{1}{k+1}\left[1-\frac{(k-1)}{k^{2m'+4}}\right] = \frac{1}{k+1} - \frac{k-1}{k+1}\cdot\frac{1}{k^{2m'+4}} > \frac{1}{k+1} - \varepsilon$$

(A.4)

Thus, we have proved condition *(i)* for both cases *a)* and *b)*. As for condition *(ii)*, in both cases *a)* and *b)* we can select the integer *p* as

$$p = \frac{k^{2m'+4}-1}{k-1}$$

(A.5)

Then (as in Eq. A.4)

$$\lim_{p\to\infty}\frac{l_{p\sigma}}{p+1} = \lim_{m'\to\infty}\frac{k^{2m'+4}-1}{k^2-1}\cdot\left(\frac{k^{2m'+4}-1}{k-1}+1\right)^{-1} = \frac{1}{k+1},$$

(A.6)

that is if *m'* is sufficiently large, we have indeed $\frac{l_{p\sigma}}{p+1} < \frac{1}{k+1} + \varepsilon.$

To make the proof work, we should select $N = \frac{k^{2m+2}-1}{k-1}$ such that for all $m' \geq m$ both of the inequalities

$$\frac{1}{k+1}\left|\frac{k^{2m'+4}-1}{k^{2m'+3}-1}\left(1+\frac{k-1}{k^{2m'+4}}\right)^{-1} - k\right| < \frac{\varepsilon}{2}$$

(A.7)

$$\frac{1}{k+1}\left|\left(1+\frac{k-1}{k^{2m'+4}}\right)^{-1} - 1\right| < \frac{\varepsilon}{2}$$

(A.8)

are satisfied. The second part of (A.1) (the case of *lim sup*) is similarly proved. ∎

*Lemma 3. Let* $\{a_i\}, A \leq a_i \leq B$ *be a bounded sequence of real numbers such that*

$$\liminf_{i\to\infty} a_i = \alpha; \quad \limsup_{i\to\infty} a_i = \beta, \tag{A.9}$$

and let $\{n_i\} = n_1, n_2, \cdots,\ n_i > 0$ be an arbitrary sequence of integers. Define a new sequence of real numbers as

$$\{b_i\} = \underbrace{a_1, a_1, \cdots, a_1}_{n_1-times}, \underbrace{a_2, \cdots, a_2}_{n_2-times}, \cdots, \underbrace{a_k, \cdots, a_k}_{n_k-times}, \cdots \tag{A.10}$$

Then (trivially) $\{b_i\}$ is also bounded, and it has the same limes inferior and limes superior as $\{a_i\}$:

$$\liminf_{i\to\infty} b_i = \alpha; \quad \limsup_{i\to\infty} b_i = \beta. \tag{A.11}$$

*Proof.* Because of (A.9), for $\forall \varepsilon > 0$ $\exists n$ such that for $\forall m \geq n$ one has

$$a_m > \alpha - \varepsilon \tag{A.12.1}$$

$$\exists k = k(m, \varepsilon) \geq m \quad such\ that \quad a_k < \alpha + \varepsilon. \tag{A.12.2}$$

Let $N := \sum_{i=1}^{n-1} n_i + 1$, and $M \geq N$ an arbitrary integer. Then for some uniquely defined $m \geq n$ one has $\sum_{i=1}^{m-1} n_i < M \leq \sum_{i=1}^{m} n_i$ and, because of (A.12.1), $b_M = a_m > \alpha - \varepsilon$.

By (A.12.1) we have $a_{m+1} > \alpha - \varepsilon$ and there exists an index $k = k(m+1, \varepsilon)$ such that $k \geq m+1$ and $a_k < \alpha + \varepsilon$. Let $K := \sum_{i=1}^{k-1} n_i + 1$. Then $K = \sum_{i=1}^{k-1} n_i + 1 \geq \sum_{i=1}^{m} n_i + 1 > M$ and $b_K = a_k < \alpha + \varepsilon$. The second part of (A.11), concerning the *limes superior*, is similarly proved. ∎